\documentclass[10pt]{amsart}

\usepackage{amsmath,amsthm,amsfonts,amscd,amssymb}
\usepackage{hyperref,wasysym}

\newtheorem{thm}{Theorem}[section]
\newtheorem{cor}[thm]{Corollary}
\newtheorem{lem}[thm]{Lemma}

\theoremstyle{definition}

\theoremstyle{remark}
\newtheorem{rem}{Remark}[section]

\begin{document}

\title{Cyclic and well-rounded lattices}
\author{Lenny Fukshansky}
\author{David Kogan}\thanks{Fukshansky was partially supported by the Simons Foundation grant \#519058}

\address{Department of Mathematics, 850 Columbia Avenue, Claremont McKenna College, Claremont, CA 91711}
\email{lenny@cmc.edu}
\address{Institute of Mathematical Sciences, Claremont Graduate University, Claremont, CA 91711}
\email{david.kogan@cgu.edu}


\begin{abstract}
We focus on two important classes of lattices, the well-rounded and the cyclic. We show that every well-rounded lattice in the plane is similar to a cyclic lattice, and use this cyclic parameterization to count planar well-rounded similarity classes defined over a fixed number field with respect to height. We then investigate cyclic properties of the irreducible root lattices in arbitrary dimensions, in particular classifying those that are simple cyclic, i.e. generated by rotation shifts of a single vector. Finally, we classify cyclic, simple cyclic and well-rounded cyclic lattices coming from rings of integers of Galois algebraic number fields.
\end{abstract}

\maketitle

\def\A{{\mathcal A}}
\def\B{{\mathcal B}}
\def\C{{\mathcal C}}
\def\D{{\mathcal D}}
\def\E{{\mathcal E}}
\def\F{{\mathcal F}}
\def\x{{\mathcal H}}
\def\I{{\mathcal I}}
\def\J{{\mathcal J}}
\def\K{{\mathcal K}}
\def\L{{\mathcal L}}
\def\M{{\mathcal M}}
\def\O{{\mathcal O}}
\def\R{{\mathcal R}}
\def\s{{\mathcal S}}
\def\V{{\mathcal V}}
\def\W{{\mathcal W}}
\def\X{{\mathcal X}}
\def\Y{{\mathcal Y}}
\def\H{{\mathcal H}}
\def\OO{{\mathcal O}}
\def\BB{{\mathbb B}}
\def\cee{{\mathbb C}}
\def\pee{{\mathbb P}}
\def\que{{\mathbb Q}}
\def\real{{\mathbb R}}
\def\zed{{\mathbb Z}}
\def\hyp{{\mathbb H}}
\def\aa{{\mathfrak a}}
\def\qbar{{\overline{\mathbb Q}}}
\def\eps{{\varepsilon}}
\def\ahat{{\hat \alpha}}
\def\bhat{{\hat \beta}}
\def\gt{{\tilde \gamma}}
\def\h{{\tfrac12}}
\def\be{{\boldsymbol e}}
\def\bei{{\boldsymbol e_i}}
\def\bff{{\boldsymbol f}}
\def\ba{{\boldsymbol a}}
\def\bb{{\boldsymbol b}}
\def\bc{{\boldsymbol c}}
\def\bm{{\boldsymbol m}}
\def\bk{{\boldsymbol k}}
\def\bi{{\boldsymbol i}}
\def\bl{{\boldsymbol l}}
\def\bq{{\boldsymbol q}}
\def\bu{{\boldsymbol u}}
\def\bt{{\boldsymbol t}}
\def\bs{{\boldsymbol s}}
\def\bv{{\boldsymbol v}}
\def\bw{{\boldsymbol w}}
\def\bx{{\boldsymbol x}}
\def\bX{{\boldsymbol X}}
\def\bz{{\boldsymbol z}}
\def\bwy{{\boldsymbol y}}
\def\bY{{\boldsymbol Y}}
\def\bL{{\boldsymbol L}}
\def\baa{{\boldsymbol\alpha}}
\def\bbb{{\boldsymbol\beta}}
\def\bet{{\boldsymbol\eta}}
\def\bxi{{\boldsymbol\xi}}
\def\bo{{\boldsymbol 0}}
\def\b1{{\boldsymbol 1}}
\def\bol{{\boldkey 1}_L}
\def\ep{\varepsilon}
\def\p{\boldsymbol\varphi}
\def\q{\boldsymbol\psi}
\def\b1{\boldsymbol 1}
\def\rank{\operatorname{rank}}
\def\aut{\operatorname{Aut}}
\def\lcm{\operatorname{lcm}}
\def\sgn{\operatorname{sgn}}
\def\spn{\operatorname{span}}
\def\md{\operatorname{mod}}
\def\Norm{\operatorname{Norm}}
\def\dim{\operatorname{dim}}
\def\det{\operatorname{det}}
\def\Vol{\operatorname{Vol}}
\def\rk{\operatorname{rk}}
\def\Gal{\operatorname{Gal}}
\def\WR{\operatorname{WR}}
\def\WO{\operatorname{WO}}
\def\GL{\operatorname{GL}}
\def\Rho{\operatorname{P}}
\def\St{\operatorname{St}}
\def\Aut{\operatorname{Aut}}
\def\Tr{\operatorname{Tr}}

\section{Introduction}
\label{intro}

Let $n \geq 2$, and write $\|\ \|$ for the Euclidean norm on $\real^n$. Let $L \subset \real^n$ be a lattice of rank $n$, and define its {\it minimum} to be
$$|L| := \min \left\{ \|\bx\| : \bx \in L \setminus \{ \bo \} \right\}.$$
Then the set of {\it minimal vectors} of $L$ is
$$S(L) := \left\{ \bx \in L : \|\bx\| = |L| \right\}.$$
Let us write $\real_+$ for the set of all positive real numbers and $\O_n(\real)$ for the group of real orthogonal $n \times n$ matrices. We define the equivalence relation of {\it similarity} on lattices in $\real^n$ as follows: two lattices $L_1$ and $L_2$ are called {\it similar}, denoted $L_1 \sim L_2$, if there exists $\alpha \in \real_+$ and $U \in \O_n(\real)$ such that $L_2 = \alpha U L_1$. This is an equivalence relation on the space of lattices in $\real^n$, and we will write $\left< L \right>$ for the similarity class of $L$. 

We say that $L$ is {\it well-rounded} (abbreviated WR) if $\spn_{\real} S(L) = \spn_{\real} L$, $L$ is {\it generated by minimal vectors} if $L = \spn_{\zed} S(L)$, and $L$ {\it has a basis of minimal vectors} if $S(L)$ contains a basis for $L$. Notice that $S(L_1)$ is taken to $S(L_2)$ under similarity, and hence these conditions are preserved. Thus we write $\WR_n$ for the set of similarity classes of WR lattices in $\real^n$, $\WR'_n$ for the set of similarity classes generated by minimal vectors, and $\WR''_n$ for the set of similarity classes having a basis of minimal vectors. Then
$$\WR''_n \subseteq \WR'_n \subseteq \WR_n,$$
where the first containment is known to be proper for all $n \geq 10$ and the second containment is proper for all $n \geq 5$. WR lattices are central objects in lattice theory and discrete optimization; see~\cite{martinet} for many more details, as well as~\cite{mcmullen}.

A lattice $L \subset \real^n$ of rank $k$ can be written as $L = B\zed^k$, where $1 \leq k \leq n$ and $B$ is an $n \times k$ basis matrix of rank~$k$. The {\it determinant} of $L$ is then defined as
$$\det(L) = \sqrt{ \det (B^{\top} B) }.$$
A lattice $L$ is called {\it semi-stable} if for any sublattice $M$ of rank $1 \leq k \leq n$,
$$\det(L)^{1/n} \leq \det(M)^{1/k},$$
and $L$ is called {\it stable} if this inequality is strict for all $M \neq L$. The stability condition is again preserved under similarity, and we write $\St_n$ for the set of semi-stable similarity classes in~$\real^n$. Semi-stable lattices, alongside well-rounded lattices, are of great interest in lattice theory and reduction theory (see~\cite{casselman} and~\cite{shapira_weiss}).

There are two more interesting classes of lattices we would like to introduce. Let $B = \{ \bb_1,\dots,\bb_n \}$ be an ordered basis for a lattice $L$ in~$\real^n$, and define a sequence of angles $\theta_1,\dots,\theta_{n-1}$ as follows: each $\theta_i$ is the angle between $\bb_{i+1}$ and the subspace $\spn_{\real} \{ \bb_1,\dots,\bb_i \}$. Then each $\theta_i \in (0,\pi/2]$. The basis $B$ is called {\it weakly nearly-orthogonal} if $\theta_i \geq \pi/3$ for each $1 \leq i \leq n-1$. A basis $B$ is called {\it nearly-orthogonal} if every ordering of it is weakly nearly-orthogonal. A lattice $L$ is called (weakly) nearly-orthogonal if it has a (weakly) nearly-orthogonal basis. This property is again preserved under similarity, and we will write $\W_n$ (respectively, $\W^*_n$) for the set of weakly nearly-orthogonal (respectively, nearly orthogonal) well-rounded similarity classes. Nearly-orthogonal lattices have applications in image compression and signal processing~\cite{baraniuk}. Finally, a lattice $L \subset \real^n$, not necessarily of full rank, is called {\it cyclic} in $\real^n$ if it closed under the rotation shift linear operator $\rho : \real^n \to \real^n$, given by 
\begin{equation}
\label{rho}
\rho(c_1,c_2,\dots,c_n) = (c_n,c_1,\dots,c_{n-1}),
\end{equation}
i.e. if $\rho(L) = L$. This property is not preserved under similarity: indeed, the integer lattice $\zed^2$ is cyclic and is similar to $\begin{pmatrix} 1 & -a \\ a & 1 \end{pmatrix} \zed^2$, which is not cyclic for any irrational $a$. On the other hand, a full-rank lattice $L \subset \real^n$ is similar to a cyclic lattice if and only if $L$ has an isometry with minimal polynomial $x^n-1$. Cyclic lattices have been especially studied in the context of lattice-based cryptography, e.g.~\cite{mic1}, \cite{peikert}.

Reduction theory aims to specify some ``canonical" choice of representatives of similarity classes, and from this point of view it is interesting to understand the relation between the classes of lattices defined above. In general, these conditions on lattices are independent: stable lattices do not need to be well rounded or nearly-orthogonal, and there are families of well-rounded nearly-orthogonal lattices that are not semi-stable (see Lemma~1.1 of~\cite{lf:stable}). Further, well-rounded lattices are also not necessarily nearly-orthogonal (as is evident from multiple constructions and examples in~\cite{lf:dk}). The situation, however, is much simpler in dimension two: here we have the following chain of inclusions:
$$\WR''_2 = \WR'_2 = \WR_2 = \W^*_2 = \W_2 \subsetneq \St_2.$$
One can ask where do cyclic lattices fit in this picture. This is our first observation.

\begin{thm} \label{cyclic_WR} Every WR lattice $L \subset \real^2$ is similar to a unique cyclic lattice
\begin{equation}
\label{Mx}
M(x) = \begin{pmatrix} 1 & x \\ x & 1 \end{pmatrix} \zed^2
\end{equation}
with $x \in [0,2-\sqrt{3}]$. Further, if $K \subseteq \real$ is a subfield such that $L \subset K^2$, then $x \in K$.
\end{thm}

\noindent
We prove Theorem~\ref{cyclic_WR} in Section~\ref{circulant}, where our main tool is the circulant preconditioner for arbitrary matrices originally defined by Tony Chan~\cite{tony_chan} in the context of certain numerical linear optimization problems. Notice, on the other hand, that the converse of Theorem~\ref{cyclic_WR} is not true: not all cyclic lattices in the plane are WR. In fact, not all of them are even stable: for instance, the cyclic lattice $\begin{pmatrix} 3 & 2 \\ 2 & 3 \end{pmatrix} \zed^2$ is not stable.

Using Theorem~\ref{cyclic_WR}, for each WR lattice $L \subset \real^2$, let $x_L \in [0,2-\sqrt{3}]$ be the unique real number so that $L \sim M(x_L)$ as in~\eqref{Mx}. This description of similarity classes allows for a way to count them. Given a real number field $K$, we will say that a lattice $L \subset \real^2$ is {\it defined over} $K$ if $L \subset K^2$. Further, we say that a similarity class is defined over $K$ if it contains a lattice defined over $K$. Theorem~\ref{cyclic_WR} guarantees that if a well-rounded lattice $L$ is defined over $K$, then so is $M(x_L)$. Then WR similarity classes defined over $K$ are precisely those containing $M(x)$ as in~\eqref{Mx} with $x \in K$. We can define the height of a similarity class $\left< L \right>$ defined over $K$ to be the Weil height of $x_L$, denoted $H(\left< L \right>)$. With this notation, we can prove the following estimate.

\begin{thm} \label{sim_count} Let $K$ be a real number field of degree $d$, then for any $T \geq 1$,
$$\left| \left\{ \left< L \right> \text{defined over } K : H(\left< L \right>) \leq T \right\} \right| \leq \frac{\pi}{2 \sqrt{12}} \left( 1+4^{2(d+1)} \left( 2+\sqrt{3} \right)^d T^{2d} \right).$$
\end{thm}

We review all the necessary notation of height functions and prove Theorem~\ref{sim_count} in Section~\ref{ht_count}. Our main tool there is a counting estimate for algebraic numbers of bounded height due to Loher and Masser~\cite{loher_masser}. Notice that introducing the height machinery allows for explicit counting: any set of points of explicitly bounded height over a fixed number field is necessarily finite by Northcott's theorem. Indeed, our approach here is different from some previous counting estimates on planar well-rounded lattices, where the lattices in question would be taken to be sublattices of a fixed lattice in the plane and counted with respect to index (see \cite{lf1}, \cite{lf2}, \cite{lf3}, \cite{lf4}, \cite{baake}, \cite{kuhnlein}). Instead, we are counting {\it all} similarity classes defined over a fixed number field. In contrast, arithmetic similarity classes of well-rounded lattices have been counted with respect to a somewhat differently defined height in~\cite{lfpf}: these lattices are defined over quadratic number fields. We compare the estimate obtained in~\cite{lfpf} with our Theorem~\ref{sim_count} in Section~\ref{ht_count}.

In higher dimensions well-rounded lattices cannot be so nicely parameterized by cyclic ones (Lemma~\ref{dim3}). This being said, there are plenty of important lattices that are cyclic. Indeed, the condition that a lattice $L$ of rank $n$ is cyclic is equivalent to the condition that $\Aut(L)$, the automorphism group of $L$, contains the permutation matrix corresponding to the standard $n$-cycle $(1 \dots n)$, which is not trivial (for a generic lattice $\Aut(L) = \{\pm I_n \}$), and lattices with large automorphism groups are of special interest in lattice theory and the arithmetic theory of quadratic forms. We distinguish a special subclass of cyclic lattices: we will say that a lattice~$L$ of rank $n$ is {\it simple cyclic} if there exists $\ba \in L$ so that
$$L = \Lambda(\ba) := \spn_{\zed} \left\{ \ba, \rho(\ba), \dots, \rho^{n-1}(\ba) \right\},$$
i.e. simple cyclic lattices are generated by the rotation shifts of a single vector. We discuss basic properties of cyclic lattices in more details in Section~\ref{cyclic}. In Section~\ref{cyclic_root} we prove the following observation on the cyclic properties of the root lattices.

\begin{thm} \label{root_thm} The following statements hold for the root lattices and their duals:
\begin{enumerate}

\item For every $n \geq 2$, the root lattice $A_n$ and its dual $A_n^*$ are both simple cyclic lattices of rank $n$ in $\real^{n+1}$.

\item For every $n \geq 2$, the root lattice $D_n$ and its dual $D_n^*$ are both cyclic of full rank in $\real^n$. Further, $D_n$ and $D_n^*$ are simple cyclic if and only if $n$ is odd.

\item The self-dual root lattice $E_8$ is cyclic, but not simple cyclic in $\real^8$.

\item The root lattices $E_6$, $E_7$ are non-cyclic sublattices of $E_8$ in $\real^8$.

\end{enumerate}
\end{thm}

\noindent
We briefly recall the definitions and necessary properties of the classical root lattices in Section~\ref{cyclic_root} before proving Theorem~\ref{root_thm} in a series of lemmas. Finally, in Section~\ref{nf_lattices} we focus on lattices coming from rings of integers of Galois number fields. We prove the following result, where the lattices in question are viewed as cyclic under the rotational shift operator as in~\eqref{rho}, but on~$\cee^n$ instead of~$\real^n$.

\begin{thm} \label{Galois} Let $K$ be a Galois number field and $\Lambda_K$ be the lattice coming from the ring of integers $\O_K$ via a standard embedding into $K_{\real} = K \otimes_{\que} \real \subseteq \cee^d$. Then $\Lambda_K$ is cyclic in~$\cee^d$ if and only if $K/\que$ is a cyclic extension. Further, it is simple cyclic if and only if $K/\que$ is tamely ramified. In particular, a cyclotomic lattice $\Lambda_{\que(\zeta_n)}$ with $\zeta_n = e^{2\pi i/n}$ is cyclic if and only if $n=2,4,p^k$ or $2p^k$ for an odd prime $p$ and integer $k \geq 1$, and it is simple cyclic if and only if $n=2, p$, or $2p$.
\end{thm}

\noindent
We recall all the necessary number field notation in Section~\ref{nf_lattices}. Further, we comment on well-roundness properties of such cyclic lattices $\Lambda_K$ and discuss some non-cyclotomic examples (Remark~\ref{more_cyclic}). We are now ready to proceed.
\bigskip

\section{Approximations by circulant matrices}
\label{circulant}

In this section we define circulant approximation to a matrix and use it to prove Theorem~\ref{cyclic_WR}. For an $n \times n$ real matrix $A$, write $\|A\|$ for its Frobenius norm, i.e. the Euclidean norm of $A$ viewed as a vector in $\real^{n^2}$. Let $\left<\ ,\ \right>$ stand for the corresponding inner product on vectors and on matrices. For each vector $\bc = (c_1,\dots,c_n) \in \real^n$, let 
\begin{equation}
\label{Pc}
\Rho(\bc) = \begin{pmatrix} c_1 & c_2 & \dots & c_n \\ c_n & c_1 & \dots & c_{n-1} \\ \vdots & \vdots & \ddots & \vdots \\ c_2 & c_3 & \dots & c_1 \end{pmatrix}
\end{equation}
be the corresponding $n \times n$ circulant matrix. Let
$$\Pi_n = \begin{pmatrix} 0 & 1 & 0 & \dots & 0 \\ 0 & 0 & 1 & \dots & 0 \\ \vdots & \vdots & \vdots & \ddots & \vdots \\ 0 & 0 & 0 & \dots & 1 \\ 1 & 0 & 0 & \dots & 0  \end{pmatrix}$$
be a permutation matrix of order $n$, then $\Pi_n^k$ for $0 \leq k \leq n-1$ is also a permutation matrix. 

We can now define a circulant approximation to an $n \times n$ matrix $A$ (also called a circulant preconditioner), as in~\cite{tony_chan}, \cite{chu}, by $\Rho(A) := \Rho(c_1,\dots,c_n)$, where
$$c_k = \frac{1}{n} \left< A, \Pi_n^{k-1} \right>,$$
for each $1 \leq k \leq n$; in other words, each entry $c_k$ of this circulant matrix is the average of the corresponding diagonal of $A$ wrapped around to extend to full length. This circulant approximation $\Rho(A)$ was introduced by T. Chan~\cite{tony_chan}, who showed that it minimizes $\|A-C\|$ among all circulant matrices $C$.

Let $L \subset \real^n$ be a full-rank lattice with a basis $\{ \ba_1,\dots,\ba_n \}$. Write $A = (\ba_1\ \dots\ \ba_n)$ for the corresponding basis matrix. Define $\Rho_A(L) = \Rho(A^{\top})^{\top} \zed^n$ to be the corresponding cyclic lattice, which we call a {\it cyclic approximation} to $L$. Let us consider this construction for $n=2$, in which case $L = A\zed^2$, where
\begin{equation}
\label{PAL}
A = (\ba_1\ \ba_2) = \begin{pmatrix} a_{11} & a_{21} \\ a_{12} & a_{22} \end{pmatrix}, \text{ and so }\Rho_A(L) = \frac{1}{2} \begin{pmatrix} a_{11} + a_{22} & a_{12} + a_{21} \\ a_{12} + a_{21} & a_{11} + a_{22} \end{pmatrix} \zed^2.
\end{equation}
We can assume that at least one of $a_{11} + a_{22}$, $a_{12} + a_{21}$ is nonzero: if both of them are, replace $\ba_1$ with $-\ba_1$ (or $\ba_2$ with $-\ba_2$).

\begin{lem} \label{WR_cyclic} Suppose that $L \subset \real^2$ is well-rounded and $A$ is a minimal basis matrix. Then $L$ is similar to $\Rho_A(L)$.
\end{lem}

\proof
Assuming $A$ as above is a minimal basis matrix is equivalent to saying that
\begin{equation}
\label{wr-cond}
\|\ba_1\|^2 = a_{11}^2 + a_{12}^2 = a_{21}^2 + a_{22}^2 = \|\ba_2\|^2,
\end{equation}
and cosine of the angle $\theta(\ba_1,\ba_2)$ between $\ba_1$ and $\ba_2$ satisfies
$$|\cos \theta(\ba_1,\ba_2)| = \frac{|\left<\ \ba_1,\ba_2 \right>|}{\|\ba_1\| \|\ba_2\|} = \frac{|a_{11}a_{21}+a_{12}a_{22}|}{a_{11}^2 + a_{12}^2} \leq \frac{1}{2}.$$
Let us write 
$$B = (\bb_1\ \bb_2) = \begin{pmatrix} a_{11} + a_{22} & a_{12} + a_{21} \\ a_{12} + a_{21} & a_{11} + a_{22} \end{pmatrix},$$
so $\Rho_A(L) = \frac{1}{2} B \zed^2$. Clearly $\|\bb_1\| = \|\bb_2\|$, and so $L$ is similar to $\Rho_A(L)$ if and only if
$$|\cos \theta(\ba_1,\ba_2)| = |\cos \theta(\bb_1,\bb_2)|.$$
Observe a basic algebraic identity for any four real numbers $p,q,s,t$:
\begin{equation}
\label{pqst}
\frac{p+q}{s+t} = \frac{p}{s} \text{ if and only if } pt = qs.
\end{equation}
Using \eqref{pqst} along with \eqref{wr-cond}, we notice that
\begin{eqnarray*}
\cos \theta(\bb_1,\bb_2) & = & \frac{2(a_{11} + a_{22})(a_{12} + a_{21})}{(a_{11} + a_{22})^2 + (a_{12} + a_{21})^2} \\
& = & \frac{(a_{11}a_{21}+a_{12}a_{22}) + (a_{11}a_{12}+a_{21}a_{22})}{(a_{11}^2+a_{12}^2) + (a_{11}a_{22}+a_{12}a_{21})} = \cos \theta(\ba_1,\ba_2),
\end{eqnarray*}
because
$$(a_{11}a_{21}+a_{12}a_{22})(a_{11}a_{22}+a_{12}a_{21}) = (a_{11}a_{12}+a_{21}a_{22})(a_{11}^2+a_{12}^2).$$
Hence $L$ is similar to $\Rho_A(L)$.
\endproof

\proof[Proof of Theorem~\ref{cyclic_WR}]
Write $L = A\zed^2$ as in~\eqref{PAL}, then by Lemma~\ref{WR_cyclic}, $L$ is similar to $\Rho_A(L)$ as in~\eqref{PAL}. Notice that 
\begin{equation}
\label{nz_cond}
(a_{11} + a_{22})(a_{12} + a_{21}) = 0
\end{equation}
if and only if $\Rho_A(L)$ is similar to $\zed^2$, which is precisely $M(x)$ as in~\eqref{Mx} with $x=0$. On the other hand, \eqref{nz_cond} does not hold if and only if $\Rho_A(L)$ is similar to $M(x)$ with 
\begin{equation}
\label{x-def}
x = \frac{a_{12} + a_{21}}{a_{11} + a_{22}}.
\end{equation}
Suppose now that $x \neq y$ are such that 
$$\begin{pmatrix} 1 & x \\ x & 1 \end{pmatrix},\ \begin{pmatrix} 1 & y \\ y & 1 \end{pmatrix}$$
are minimal basis matrices for $M(x)$ and $M(y)$, respectively, and $M(x)$ is similar to $M(y)$. Then absolute values of cosines of the angles between these minimal basis vectors are equal and $\leq 1/2$, in particular
$$\frac{2|x|}{x^2+1} = \frac{2|y|}{y^2+1}.$$
This is true if and only if $y = -x$ or $y = \pm 1/x$. Hence, to ensure uniqueness, we can take $0 < x \leq 1$. Additionally, we need
$$\frac{2x}{x^2+1} \leq \frac{1}{2},$$
which means $0 < x \leq 2-\sqrt{3}$. Combining this with the case $x=0$ completes the proof of the first part of the theorem. The second part follows from the fact that $x$ is given by either $\pm$ the expression in~\eqref{x-def} or its inverse, which therefore lies in the same field that contains $a_{ij}$'s.
\endproof

\begin{rem} \label{WR_par} There is a standard parameterization of well-rounded similarity classes in the plane by lattices of the form
$$\Lambda(a,b) = \begin{pmatrix} 1 & a \\ 0 & b \end{pmatrix} \zed^2,$$
where $0 \leq a \leq 1/2$ with $a^2+b^2=1$ (see, for instance,~\cite{lfpf}). However, if a lattice $L$ is similar to some such $\Lambda(a,b)$, they are not necessarily defined over the same field, unlike the parameterization of our Theorem~\ref{cyclic_WR}.
\end{rem}

Our observations above imply, in particular, that in $\real^2$ well-rounded lattices are always similar to cyclic lattices. Unsurprisingly, this is not true in higher dimensions.

\begin{lem} \label{dim3} A well-rounded lattice of rank $\geq 3$ is not necessarily similar to a cyclic lattice.
\end{lem}

\proof
We can give a simple dimensional argument. First, consider the set $\WR''_n$, i.e. similarity classes of WR lattices with a basis of minimal vectors. For all $n \geq 2$, the set $\WR''_n$ is determined by $n \times n$ matrices $A = (a_{ij})$ with 
$$a_{11} = 1,\ a_{21} = \dots = a_{n1} = 0,\ \sum_{j=1}^n a_{ij}^2 = 1\ \forall\ 2 \leq i \leq n,$$
also satisfying some inequalities. Since inequalities do not reduce the dimension, this space has dimension 
$$n^2 - n - (n-1) = (n-1)^2.$$
Of course, for $n \geq 5$, $\WR''_n \neq \WR_n$, however each WR lattice $L$ has only finitely many sublattices with a basis consisting of some minimal vectors of $L$, and each WR lattice with a basis of minimal vectors is contained in only finitely many WR lattices with these vectors among the minimal (see~\cite{martinet_mink}). Hence, the sets $\WR_n$ and $\WR''_n$ have the same dimension.

On the other hand, the space of similarity classes of cyclic lattices in $\real^n$ is determined by $n \times n$ matrices $A = (a_{ij})$ with 
$$a_{11} = 1,\ a_{ij} = \rho^j(a_{i1})\ \forall\ 2 \leq i \leq n.$$
This space has dimension
$$n-1 < (n-1)^2\ \forall\ n \geq 3.$$
Thus the space of WR similarity classes is too big to be parameterized by cyclic lattices. Notice, however, that when $n=2$ these dimensions coincide.
\endproof

\bigskip

\section{Counting WR similarity classes}
\label{ht_count}

Let $K$ be a number field of degree $d := [K : \que] \geq 1$. We write
$$M(K) = M_{\infty}(K) \cup M_f(K),$$
for the set of places of $K$, split into the subsets $M_{\infty}(K)$ of archimedean and $M_f(K)$ of non-archimedean places. The archimedean places correspond to the embeddings
$$\sigma_1,\dots,\sigma_d : K \hookrightarrow \cee$$
of $K$ as usual: $M_{\infty}(K) = \{ v : v = v(\sigma_i) \text{ for some } 1 \leq i \leq d\}$, where for each $1 \leq i \leq d$ and $x \in K$,
$$|x|_{v(\sigma_i)} = |\sigma_i(x)| = |\bar{\sigma}_i(x)| = |x|_{v(\bar{\sigma}_i)}$$
since complex conjugate embeddings give rise to the same place (we regard places in $M_{\infty}(K)$ without repetition). We order the embeddings so that $\sigma_1$ extends to the identity map on $\cee$ and so $v_1 = v(\sigma_1)$ is the place corresponding to it. For each $v \in M(K)$ let $d_v = [K_v : \que_v]$ be the local degree, then for each $u \in M(\que)$, $\sum_{v \mid u} d_v = d$. We normalize the absolute values so that for each nonzero $x \in K$ the product formula reads
$$\prod_{v \in M(K)} |x|_v^{d_v} = 1.$$
Let $n \geq 2$, and for any place $v \in M(K)$ and $\bx = (x_1,\dots,x_n) \in K^n$ we define the corresponding sup-norm $|\bx|_v = \max \{ |x_1|_v,\dots,|x_n|_v \}$. Then the height $H : K^n \to \real_{\geq 0}$ is given by
$$H(\bx) = \prod_{v \in M(K)} |\bx|_v^{\frac{d_v}{d}}.$$
The Weil height $h : K \to \real_{\geq 1}$ is then defined by $h(x) = H(1,x)$. This height is absolute, meaning that $H(\bx)$ is the same when computed over any number field $K$ containing the coordinates of $\bx$: this is due to the normalizing exponent~$1/d$ in the definition.

We also define local and ``anti-local" heights following~\cite{loher_masser}. Let $\bx \in \qbar^n$ and let $K'$ be an extension of $K$ containing the coordinates of $\bx$. For each archimedean $v \in M(K)$, define
$$H_v(\bx) = \prod_{w \in M(K'), w \mid v} |\bx|_w^{\frac{d_w}{[K':\que]}},\ H^v(\bx) = \prod_{w \in M(K'), w \nmid v} |\bx|_w^{\frac{d_w}{[K':\que]}},$$
so that $H(\bx) = H_v(\bx) H^v(\bx)$. 

\begin{lem} \label{ht_bnd} Let $K$ be a number field of degree $d$ and $\alpha \in \qbar$. Let $K' = K(\alpha)$ and let $u \in M_{\infty}(K)$. For real $T \geq 1$, define
$$S_K(\alpha,T) = \left\{ x \in K : |x|_u \leq |\alpha|_u, h(x) \leq T \right\}.$$
Then
$$|S_K(\alpha,T)| \leq \frac{\pi}{\sqrt{12}} \left( 1+4^{2(d+1)} (Th(\alpha))^{2d} \right).$$
\end{lem}

\proof
Define
$$S'_K(\alpha,T) = \left\{ x \in K : |x|_u \leq |\alpha|_u, H(\alpha,x) \leq T \right\},$$
and notice that for $x \in S_K(\alpha,T)$, we have
$$T h(\alpha) \geq h(x) h(\alpha) \geq H(\alpha,x),$$
which means that $|S_K(\alpha,T)| \leq |S'_K(\alpha,T h(\alpha))|$. Notice that
$$H(\alpha,x) = H_u(\alpha,x) H^u(\alpha,x),$$
and so
\begin{eqnarray*}
|x|_u \leq |\alpha|_u \Leftrightarrow H_u(\alpha,x) & = & \prod_{w \in M(K'), w \mid u} \max \{ |\alpha|_w, |x|_w \}^{\frac{d_w}{[K':\que]}} \\
& \leq & |\alpha|^{\frac{1}{[K':\que]} \sum_{w \mid u} d_w}_u = |\alpha|^{\frac{[K':K]}{[K':\que]}}_u = H_u(\alpha).
\end{eqnarray*}
Hence
$$S'_K(\alpha,Th(\alpha)) = \left\{ x \in K : H_u(\alpha,x) \leq H_u(\alpha), H(\alpha,x) \leq Th(\alpha) \right\},$$
and so if
\begin{equation}
\label{lm_cond}
H_u(\alpha,x) \leq H_u(\alpha),\ H^u(\alpha,x) \leq Th(\alpha) H^u(\alpha),
\end{equation}
then $x \in S'_K(\alpha,Th(\alpha))$, since $H(\alpha) = H^u(\alpha) H^u(\alpha) = 1$ by the product formula. Let $N_K(\alpha,T)$ be the number of elements $x \in K$ satisfying conditions~\eqref{lm_cond}, then
$$|S_K(\alpha,T)| \leq |S'_K(\alpha,T h(\alpha))| \leq N_K(\alpha,T),$$
and by the Proposition in Section~3 of~\cite{loher_masser},
$$N_K(\alpha,T) \leq \frac{\pi}{\sqrt{12}} \left( 1+4^{2(d+1)} (Th(\alpha))^{2d} \right).$$
This completes the proof of the lemma.
\endproof

\proof[Proof of Theorem~\ref{sim_count}]
By Theorem~\ref{cyclic_WR},
$$\left| \left\{ \left< L \right> \text{defined over } K : H(\left< L \right>) \leq T \right\} \right| = \frac{1}{2} \left| \left\{ x \in K : |x|_{v_1} \leq 2-\sqrt{3}, h(x) \leq T \right\} \right|,$$
where $1/2$ accounts for the fact that we are only considering positive $x$. The theorem then follows from Lemma~\ref{ht_bnd} combined with the fact that
\begin{eqnarray*}
h(2-\sqrt{3}) & = & \prod_{v \in M(\que(\sqrt{3}))} \max \left\{ 1, |2-\sqrt{3}|_v \right\}^{d_v/2} = \prod_{v \mid \infty} \max \left\{ 1, |2-\sqrt{3}|_v \right\}^{1/2} \\
& = & \max \left\{ 1, 2-\sqrt{3} \right\}^{1/2} \max \left\{ 1, 2+\sqrt{3} \right\}^{1/2} = \sqrt{2+\sqrt{3}}.
\end{eqnarray*}
\endproof

\begin{rem} \label{ht_compare} We can compare the estimate of Theorem~\ref{sim_count} in case $K$ is a quadratic number field to the estimate obtained in~\cite{lfpf}. As indicated in~\cite{lfpf}, an arithmetic well-rounded lattice in the plane is of the form
$$L(a,b) = \begin{pmatrix} 1 & \frac{a}{b} \\ 0 & \frac{\sqrt{b^2-a^2}}{b} \end{pmatrix} \zed^2,$$
where $a,b \in \zed$ are relatively prime with $0 < a \leq b/2$, or $a=0, b=1$. These are precisely WR lattices similar to those with integer-valued quadratic norm forms. Such a lattice $L(a,b)$ is similar to the cyclic lattice $M(x)$ as in~\eqref{Mx} with
\begin{equation}
\label{xab}
x = \frac{a}{\sqrt{b^2-a^2}+b}.
\end{equation}
Notice that for $x$ as in~\eqref{xab}, $h(x)$ is of the order of magnitude $O(\sqrt{b})$. Then bounding $b$ by $T$, the upper bound of our Theorem~\ref{sim_count} grows like $O(T^2)$, which is consistent with the growth order of $N_3(T)$ in Theorem~1.1 of~\cite{lfpf}. On the other hand, $N_3(T)$ counts {\it all} arithmetic lattices with $b \leq T$, whereas our Theorem~\ref{sim_count} counts only those defined over a fixed number field. The difference, however, is that the set of WR lattices defined over a fixed quadratic field is far more general than those corresponding to $x$ as in~\eqref{xab}, i.e. not nearly all of them are arithmetic. For instance, WR lattices
$$\begin{pmatrix} 1 & \frac{1}{\sqrt{5}} \\ 0 & \frac{2}{\sqrt{5}} \end{pmatrix} \zed^2 \sim \begin{pmatrix} 1 & \frac{1}{2+\sqrt{5}} \\  \frac{1}{2+\sqrt{5}} & 1 \end{pmatrix} \zed^2$$
are defined over $\que(\sqrt{5})$, but are not arithmetic. 
\end{rem}

\bigskip

\section{Cyclic lattices}
\label{cyclic}

In this section we discuss some basic properties of cyclic lattices. Let $\bc \in \real^n$ be a nonzero vector and $\Rho(\bc)$ be the corresponding circulant matrix as in~\eqref{Pc}. Let us write
\begin{equation}
\label{cx}
\bc(x) = \sum_{k=1}^n c_k x^{k-1}
\end{equation}
for the polynomial of degree $n-1$ with $\bc$ for the coefficient vector. It is a well-known fact that
\begin{equation}
\label{detc}
\det \Rho(\bc) = \prod_{j=1}^n \bc(\omega_n^j),
\end{equation}
where $\omega_n$ is a primitive $n$-th root of unity. Therefore $\Rho(\bc)$ is singular if and only if $\bc(x)$ is divisible by some cyclotomic polynomial $\Phi_d(x)$ for $d \mid n$. Otherwise, the simple cyclic lattice $\Lambda(\bc) = \Rho(\bc)^{\top} \zed^n$ has full rank. 

\begin{lem} \label{dual} A full-rank lattice $L \subset \real^n$ is cyclic if and only if its dual $L^*$ is cyclic. Further, $L$ is simple cyclic if and only if $L^*$ is simple cyclic.
\end{lem}

\proof
Recall that
$$L^* = \left\{ \bx \in \real^n : \left< \bx, \bwy \right> \in \zed\ \forall\ \bwy \in L \right\}.$$
Assume that $L$ is cyclic. Let $\bx \in L^*$ and $\bwy \in L$, then
$$\left< \rho(\bx), \bwy \right> = x_n y_1 + x_1 y_2 + x_2 y_3 + \dots + x_{n-1} y_n = \left< \bx, \rho^{-1}(\bwy) \right> \in \zed,$$
since $\rho^{-1}(\bwy) = \rho^{n-1}(\bwy) \in L$. Thus $\rho(\bx) \in L^*$ for every $\bx \in L^*$, hence $L^*$ is cyclic. Conversely, suppose $L^*$ is cyclic, then $L = (L^*)^*$ is cyclic by the argument above.

Now suppose $L$ is simple cyclic. Then $L = \Rho(\bc)^{\top} \zed^n$ for some $\bc \in L$, and so
$$L^* = \left( \Rho(\bc)^{\top}\right)^{-\top} \zed^n = \Rho(\bc)^{-1} \zed^n.$$
Since the transpose and inverse of a circulant matrix are both also circulant, we can conclude that $L^*$ has a circulant basis matrix, and hence $L^*$ is also simple cyclic. Conversely, if $L^*$ is simple cyclic, then $L = (L^*)^*$ is simple cyclic by the argument above.
\endproof

Cyclic sublattices of $\zed^n$ can be constructed algebraically. Indeed, let $R_n$ be the quotient ring $\zed[x]/ \left< x^n-1 \right>$ and define a map $\phi : R_n \to \zed^n$ which sends a polynomial $\bc(x) = \sum_{k=1}^n c_k x^{k-1} \in R_n$ to its vector of coefficients $\bc \in \zed^n$. The map $\phi$ is a free $\zed$-module isomorphism, which maps ideals in $R_n$ to sublattices in $\zed^n$. In fact, a sublattice $L \subseteq \zed^n$ is cyclic if and only if $L = \phi(I)$ for some ideal $I \subseteq R_n$: the cyclic rotation operator $\rho$ on $\zed^n$ corresponds to multiplication by $x$ in $R_n$, i.e.,
$$\phi(x \bc(x)) = \rho(\bc).$$
Further details on cyclic sublattices of $\zed^n$ that were extensively studied in the context of lattice cryptography can be found in~\cite{mic1} and~\cite{peikert}. In fact, cyclic sublattices of $\zed^n$ are a special case of the more general class of ideal lattices from quotient polynomial rings (see~\cite{lub_mic} and~\cite{ferrari} for more details on these). The following simple observation will be useful to us (see also Propositions~2.1 and~2.2 of~\cite{ferrari}). We provide a proof here for self-containment.

\begin{lem} \label{zero_div} Let $I = \left< \bc(x) \right>$ be an ideal in $R_n = \zed[x]/ \left< x^n-1 \right>$. Then $\phi(I) = \Lambda(\bc) \subseteq \zed^n$ has full rank if and only if $\bc(x)$ is not a zero-divisor in $R$.
\end{lem}

\proof
The polynomial $\bc(x)$ is not a zero-divisor in $R_n$ if and only if for any nonzero $\ba(x) = \sum_{k=0}^n a_k x^{k-1} \in R_n$,
$$\ba(x) \bc(x) = \sum_{k=1}^n a_k x^{k-1} \bc(x) \neq 0.$$
This is equivalent to the statement that
$$\sum_{k=1}^n a_k \rho^{k-1}(\bc) \neq \bo$$
in $\Lambda(\bc) = \phi(I)$, i.e. $\bc,\rho(\bc),\dots,\rho^{n-1}(\bc)$ are linearly independent, meaning that $\Lambda(\bc)$ has full rank.
\endproof

\begin{lem} \label{simple_Z} A full-rank sublattice $L \subseteq \zed^n$ is simple cyclic if and only if $L = \phi(I)$ where $I \subseteq R_n$ is a principal ideal such that the quotient ring $R_n/I$ is finite.
\end{lem}

\proof
A sublattice $L \subseteq \zed^n$ is cyclic if and only if $L = \phi(I)$ for some ideal $I \subseteq R_n$ and
$$\left| R_n / I \right| = \left| \phi(R_n) / \phi(I) \right| = \left| \zed^n / L \right|.$$
Hence $R_n/I$ is finite if and only if $L$ is a full-rank sublattice of $\zed^n$. Further, $L = \Lambda(\bc)$ for some $\bc \in \zed^n$ if and only if $I = \left< \bc(x) \right>$, i.e. $I$ is principal.
\endproof

\begin{cor} A full-rank sublattice $L \subseteq \zed^n$ is simple cyclic if and only if $L = \phi(\left< \bc(x) \right>)$ where $\bc(x) \in R_n$ is not a zero-divisor.
\end{cor}

\proof
Combine Lemmas~\ref{zero_div} and~\ref{simple_Z}.
\endproof

\begin{rem} \label{cnt_simple} Recall that
$$x^n-1 = \prod_{d \mid n} \Phi_d(x),$$
where $\Phi_d(x)$ stands for $d$-th cyclotomic polynomial. Then it follows from~\cite{rossmann} (Theorem~A combined with Theorem~1.2) that the number of ideals in $R_n$ (hence the number of cyclic sublattices of $\zed^n$) with index $\leq T$ grows like $O \left(T (\log T)^{\tau(n)-1} \right)$ as $T \to \infty$, where $\tau(n)$ is the number of divisors of $n$ (see also Theorem~2.3 of~\cite{lf:kuhnlein} for a convenient formulation). Simple cyclic sublattices of full rank correspond to principal ideals of finite index, and hence their number also grows like $O \left(T (\log T)^{\tau(n)-1} \right)$. This follows from the proof of Theorem~2.3 of~\cite{lf:kuhnlein}, since every ideal class in $R_n$ contributes equally to the total number of ideals of bounded index.\footnote{This observation is due to Stefan K\"uhnlein.}
\end{rem}

\bigskip

\section{Cyclic representation of root lattices}
\label{cyclic_root}

In this section we focus on the cyclic properties of the standard root lattices, in particular proving Theorem~\ref{root_thm}.

\begin{lem} \label{An} The root lattice
$$A_n = \left\{ \bx \in \zed^{n+1} : \sum_{i=1}^{n+1} x_i = 0 \right\}$$
and its dual $A_n^*$ are simple cyclic lattices of rank $n$ in $\real^{n+1}$ for each $n \geq 2$.
\end{lem}

\proof
Write $\sum(\bx)$ for the sum of all the coordinates of the vector $\bx$, and notice that $\sum(\bx) = \sum(\rho(\bx))$ for any vector $\bx \in \real^n$ for any $n \geq 1$. Then the lattice $A_n$ is closed under $\rho : \real^{n+1} \to \real^{n+1}$, and hence is a cyclic lattice of rank $n$ in $\real^{n+1}$. In fact,
$$A_n = \spn_{\zed} \left\{ \ba,\rho(\ba),\dots,\rho^{n-1}(\ba) \right\}$$
for the vector $\ba = (1, -1, 0,\dots,0)^{\top}$, hence it is simple cyclic. 

Let us write $\boldsymbol 1_{n+1}$ for the vector in $\real^{n+1}$ with all the coordinates equal to $1$ and $\boldsymbol 1_{n+1}^{\perp}$ for the co-dimension one subspace of $\real^{n+1}$ orthogonal to $\boldsymbol 1_{n+1}$. Then $A_n = \boldsymbol 1_{n+1}^{\perp} \cap \zed^{n+1}$ and $A_n^*$ is the orthogonal projection of $\zed^{n+1}$ onto $\boldsymbol 1_{n+1}^{\perp}$. As described in Proposition~4.2.3 of~\cite{martinet}, $A_n^*$ is generated by the vectors
$$\bwy_i = \frac{1}{n+1} \left( (n+1) \be_i - \boldsymbol 1_{n+1} \right),\ 1 \leq i \leq n+1,$$
which are rotation shifts of $\bwy_1$. Hence $A_n^*$ is also simple cyclic in $\real^{n+1}$.
\endproof

\begin{lem} \label{Dn} The root lattice
$$D_n = \left\{ \bx \in \zed^n : \sum_{i=1}^n x_i \equiv 0\ (\md 2) \right\}$$
is cyclic in $\real^n$ for each $n \geq 2$. It is simple cyclic if and only if $n$ is odd.
\end{lem}

\proof
Since $\sum(\bx) = \sum(\rho(\bx))$ for any $\bx \in \real^n$, we see that $D_n$ is cyclic. Let $n \geq 3$ be odd, and take $\bc = (1, 1, 0, \dots, 0)^{\top} \in D_n$. We will show that $D_n = \Rho(\bc)^{\top} \zed^n$. Recall that $\det D_n = 2$, hence it is sufficient to show that $\det \Rho(\bc)^{\top} = \pm 2$. Notice that
$$\Rho(\bc)^{\top} = \begin{pmatrix} 1 & 0 & 0 & \dots & 0 & 1 \\ 1 & 1 & 0 & \dots & 0 & 0 \\ 0 & 1 & 1 & \dots & 0 & 0 \\ \vdots & \vdots & \vdots & \ddots & \vdots & \vdots \\ 0 & 0 & 0 & \dots & 1 & 0 \\ 0 & 0 & 0 & \dots & 1 & 1 \end{pmatrix}.$$
Performing the Laplace expansion along the first row and keeping in mind that $n$ is odd, we see that $\det \Rho(\bc)^{\top} = 2$, and hence $D_n = \Lambda(\bc)$.

Now suppose $n \geq 2$ is even. Arguing toward a contradiction, suppose that there exists some $\bc \in D_n$ such that $D_n = \Lambda(\bc)$. Let $\bc(x)$ be as in~\eqref{cx}, then by~\eqref{detc},
$$\det \Rho(\bc) = \prod_{j=1}^n \bc(\omega_n^j) = \pm \det D_n = \pm 2.$$
In particular, $1 = \omega_n^n$ and $-1 = \omega_n^{n/2}$ are both $n$-th roots of unity, and so $\bc(1) \bc(-1)$ is a nonzero integer. Then the remaining product
$$\prod_{j=1,\ j \neq \frac{n}{2}}^{n-1} \bc(\omega_n^j) = \pm \frac{2}{\bc(1) \bc(-1)} \in \que,$$
but this product is also an algebraic integer, hence it is in $\zed$. This means that $|\bc(1) \bc(-1)| \leq 2$. On the other hand, $\bc(1) = \sum_{i=1}^n c_i$ is even, since $\bc \in D_n$, thus $\bc(1) = \pm 2$ and so $\bc(-1) = \pm 1$. Let $\alpha$ be the sum of coefficients of $\bc(x)$ in front of even powers of $x$ and $\beta$ be the sum of coefficients in front of odd powers of $x$, then
$$2 \mid \bc(1) = \alpha + \beta,\ 2 \nmid \bc(-1) = \alpha - \beta,$$
implying that $2 \nmid \bc(1) + \bc(-1) = 2\alpha$. This is a contradiction, and hence $D_n$ is not simple cyclic.
\endproof

\begin{rem} \label{Z_Dn} The cyclic lattice $A_n$ can be constructed from the ideal $\left< x-1 \right>$ of rank $n$ in $R_{n+1} = \zed[x]/ \left< x^{n+1}-1 \right>$ and $D_n$ from an ideal of full rank in $R_n = \zed[x]/ \left< x^n-1 \right>$, as we discussed in Section~\ref{cyclic}. The lattices $A_n$ (in $\real^{n+1}$) and $D_n$ for odd $n$ are simple cyclic; the latter one can also be easily obtained from the ideal $\left< x+1 \right>$ in $R_n$ (see also Propositions~4.5 and~4.2 of~\cite{ferrari}). However, $D_n$ for even $n$ is of full rank but not simple cyclic, and hence cannot come from a principal ideal in $R_n$, by Lemma~\ref{simple_Z}. In fact, it can easily be seen as the image under $\phi$ of the ideal
$$I = \left< 2, x+1 \right> = \left< x^{n-1}+x^{n-2}, -x^{n-1}+1 \right> \subset R_n.$$
This can be compared to Proposition~4.4 of~\cite{ferrari}, where $D_n$ for even $n \geq 4$ is obtained as the image of the principal ideal $\left< x+1 \right>$ in $\zed[x]/\left< x^n+1 \right>$ under the same kind of coefficient embedding into~$\zed^n$.
\end{rem}

\begin{lem} \label{E8} The lattice $E_8$ is cyclic, but not simple cyclic. The lattices $E_6$ and $E_7$ are not cyclic.
\end{lem}

\proof
Recall that the lattice $E_8$ can be defined as
$$E_8 = D_8 \cup \left( \frac{1}{2} \b1_8 + D_8 \right).$$
Notice that the vector $\b1_8$ is invariant under $\rho$, and hence $\frac{1}{2} \b1_8 + D_8$ is closed under $\rho$, as is $D_8$. This means that $E_8$ is cyclic. On the other hand, if $\bc \in E_8$, then either $\bc \in D_8$ or $\bc \in \frac{1}{2} \b1_8 + D_8$. If $\bc \in D_8$, then $\Lambda(\bc) \subseteq D_8$, so $\Lambda(\bc) \neq E_8$. We now argue similarly to our proof of Lemma~\ref{Dn} above. Suppose that $\bc \in \frac{1}{2} \b1_8 + D_8$ is such that $\Lambda(\bc) = E_8$, then $\bc = \frac{1}{2} \b1_8 + \bc'$ for some $\bc' \in D_8$, and
\begin{equation}
\label{prod_E8}
\prod_{j=1}^8 \bc(\omega_8^j) = \pm \det E_8 = \pm 1.
\end{equation}
Notice that
$$\bc(x) = \frac{1}{2} \sum_{k=0}^7 x^k + \bc'(x) = \frac{x^8-1}{2(x-1)} + \bc'(x),$$
for all $x \neq 1$, and so $\bc(\omega_8^j) = \bc'(\omega_8^j) \in \zed$ for every $1 \leq j \leq 7$ and $\bc(1) = 4 +\bc'(1)$. In order for~\eqref{prod_E8} to hold, we must in particular have
$$\bc(1) \bc(-1) = (4+\bc'(1)) \bc'(-1) = \pm 1.$$
However, $\bc'(1)$ and $\bc'(-1)$ must both be even, since $\bc' \in D_8$. This is a contradiction, and hence $E_8$ is not simple cyclic.

The root lattices $E_7$ and $E_6$ can be described as sublattices of $E_8$ orthogonal to the vector $\be_7+\be_8$ and to the pair of vectors $\be_7+\be_8$, $\be_6+\be_8$, respectively. These lattices are not cyclic, since the spaces
$$\spn_{\real} \{ \be_7+\be_8 \},\ \spn_{\real} \{ \be_7+\be_8, \be_6+\be_8 \}$$
are not closed under $\rho$. For instance, $\bx = \be_4+\be_5 \in E_7 \cap E_6$, however $\rho(\bx) \not\in E_6$ and $\rho^2(\bx) \not\in E_7$.
\endproof

\proof[Proof of Theorem~\ref{root_thm}] The theorem follows by combining Lemmas~\ref{An}, \ref{Dn} and~\ref{E8} with Lemma~\ref{dual}.
\endproof

\bigskip

\section{Number field lattices}
\label{nf_lattices}

Yet another important class of lattices comes from rings of integers of number fields. In this section we classify those of them that are cyclic, proving Theorem~\ref{Galois}. As in Section~\ref{ht_count}, let $K$ be a number field of degree $d = r_1+2r_2$ with embeddings
$$\sigma_1,\dots,\sigma_d : K \hookrightarrow \cee,$$
where $r_1$ of them are real and $2r_2$ are complex, split into conjugate pairs. Then $K_{\real} = K \otimes_{\que} \real$ can be viewed as a subspace of $\real^{r_1} \times \cee^{2r_2} \subseteq \cee^d$, given by (up to a permutation of the coordinates)
$$\left\{ (\bx,\bwy) \in \real^{r_1} \times \cee^{2r_2} : y_{r_2+j} = \bar{y}_j\ \forall\ 1 \leq j \leq r_2 \right\} \cong \real^{r_1} \times \cee^{r_2} \subseteq \cee^d.$$
Notice that in this last containment, we identify each copy of $\real$ with the real part of the corresponding copy of $\cee$.
It is a Euclidean space with respect to the bilinear form induced by the trace-form $\left< \alpha, \beta \right>$ on the number field $K$:
\begin{equation}
\label{TR}
\left< \alpha, \beta \right> := \Tr_K(\alpha \bar{\beta}) \in \real
\end{equation}
for any $\alpha, \beta \in K$, where $\Tr_K$ stands for the usual trace map on the number field~$K$. We can define the embedding
$$\Sigma_K = (\sigma_1,\dots,\sigma_d) : K \hookrightarrow K_{\real}$$
of $K$ into $K_{\real}$. The ring of integers $\O_K$ becomes a lattice of full rank in $K_{\real}$ under this embedding, and we write $\Lambda_K$ for the image $\Sigma_K(\O_K)$. An equivalent description of $\Lambda_K$ is as a free $\zed$-module $\O_K$ equipped with the bilinear form $\left< \cdot, \cdot \right>$ we defined in~\eqref{TR}. It is easy to verify that $\left< \alpha, \beta \right>$ is equal to the usual dot product of the vectors $\Sigma_K(\alpha)$ and $\Sigma_K(\beta)$ in $K_{\real}$. We write $\Aut(\Lambda_K)$ for the automorphism group of the lattice $\Lambda_K$, i.e. the group of isometries of this trace-induced bilinear form.

\begin{lem} \label{aut} Suppose $K/\que$ is a Galois extension with the Galois group $G$. Then $G \leq \Aut(\Lambda_K)$.
\end{lem}

\proof
Notice that all the embeddings of $K$ are precisely the elements of $G$, and for every $\tau \in G$ and $\alpha,\beta \in \O_K$,
\begin{eqnarray*}
\left< \tau(\alpha), \tau(\beta) \right> & = & \Tr_K(\tau(\alpha \bar{\beta})) =  \sum_{\sigma \in G} \sigma \tau(\alpha \bar{\beta}) \\
& = & \sum_{\sigma \in G} \sigma(\alpha \bar{\beta}) = \Tr_K(\alpha \bar{\beta}) = \left< \alpha, \beta \right>,
\end{eqnarray*}
since right-multiplication by $\tau$ simply permutes elements of $G$. Therefore $G$ is a subgroup of $\Aut(\Lambda_K)$.
\endproof

Notice, however, that $\Aut(\Lambda_K)$ can be quite a bit larger than the Galois group of $K/\que$. For example, $\Lambda_{\que(i)}$ is similar to $\zed^2$, which has automorphism group of order~$8$, and $\Lambda_{\que(\sqrt{-3})}$ is similar to the hexagonal lattice, which has automorphism group of order~$12$; in both cases, Galois groups of the quadratic fields have order~$2$. Hence there are often automorphisms of the lattice that do not come from the Galois action. This observation raises a question: if $\Lambda_K$ is cyclic, does the cyclic shift operator $\rho$ necessarily come from the Galois action? In the next lemma we answer this question in the affirmative. To avoid ambiguity, in this section we view lattices as cyclic under the rotational shift operator~$\rho$ as in~\eqref{rho} but on~$\cee^d$.

\begin{lem} \label{nf_cyclic} Suppose $K/\que$ is a Galois extension with Galois group $G$. Then $\Lambda_K$ is cyclic (for an appropriate ordering of the embeddings) if and only if $K/\que$ is a cyclic extension with $G = \left< \sigma \right>$, where the automorphism $\sigma : K \to K$ is such that
\begin{equation}
\label{rho_sigma}
\rho \left( \Sigma_K(\alpha) \right) = \Sigma_K(\sigma(\alpha)),
\end{equation}
for every $\alpha \in \O_K$. 
\end{lem}

\proof
Suppose first that $\Lambda_K$ is cyclic, then for any $\alpha \in \O_K$,
$$\rho \left( \sigma_1(\alpha),\dots,\sigma_d(\alpha) \right) = \left( \sigma_d(\alpha), \sigma_1(\alpha),\dots,\sigma_{d-1}(\alpha) \right) \in \Lambda_K.$$
This means that $\sigma_d(\alpha) \in \O_K$ and
$$\sigma_1 \sigma_d(\alpha) = \sigma_d(\alpha), \sigma_2 \sigma_d(\alpha) = \sigma_1(\alpha), \dots, \sigma^2_d(\alpha) = \sigma_{d-1}(\alpha),$$
for all $\alpha \in \O_K$. Then $\sigma_1$ is the identity map, and
$$\sigma_2 \sigma_d = \sigma_1, \sigma_3 \sigma_d = \sigma_2, \sigma_4 \sigma_d = \sigma_3, \dots, \sigma^2_d = \sigma_{d-1}.$$
This implies that
$$\sigma_j = \sigma_d^{d-j+1},\ \forall\ 1 \leq j \leq d-1,$$
and the action of $\rho$ on $\Lambda_K$ is given by the action of $\sigma_d$ on $\O_K$ as specified in~\eqref{rho_sigma}.

On the other hand, suppose that $K/\que$ is cyclic with Galois group $G = \left< \sigma \right>$ so that the embeddings are ordered as
\begin{equation}
\label{order}
\sigma_j =\sigma^{d-j+1},\ \forall\ 1 \leq j \leq d-1.
\end{equation}
Then for any $\Sigma(\alpha) \in \Lambda_K$, we have
\begin{eqnarray*}
\rho \left( \alpha, \sigma^{d-1}(\alpha), \sigma^{d-2}(\alpha),\dots, \sigma(\alpha) \right) & = & \left( \sigma(\alpha), \alpha, \sigma^{d-1}(\alpha),\dots, \sigma^2(\alpha) \right) \\
& = & \Sigma_K(\sigma(\alpha)) \in \Lambda_K,
\end{eqnarray*}
and so $\Lambda_K$ is closed under $\rho$, hence is cyclic.
\endproof

In fact, one can also ask which of these cyclic lattices of the form $\Lambda_K$ are simple cyclic. We discuss this next. A normal basis for a Galois number field $K$ is a basis consisting of all conjugates of one algebraic number, and a normal integral basis for $K$ is a $\zed$-basis like this for $\O_K$. The normal basis theorem guarantees that every number field has a normal basis, however having a normal integral basis is a much more delicate property. The finite Galois extension $K/\que$ is called tamely ramified if all the ramification indices for every rational prime~$p$ are relatively prime with~$p$. The Hilbert-Speiser theorem asserts that an abelian number field (i.e. Galois number field with abelian Galois group) has a normal integral basis if and only if it is tamely ramified (see, for instance, Chapter~9 of~\cite{long} for the details).

\begin{lem} \label{K_sc} Let $K$ be a cyclic Galois number field. Then the lattice $\Lambda_K$ is simple cyclic (for an appropriate ordering of the embeddings) if and only if $K/\que$ is tamely ramified.
\end{lem}

\proof
Let $d = [K:\que]$ and $G = \left< \sigma \right>$ be the Galois group of $K/\que$. Let $\sigma_1,\dots,\sigma_d$ be the embeddings of $K$, ordered as in~\eqref{order}. By the Hilbert-Speiser theorem $K/\que$ is tamely ramified if and only if $K$ has a normal integral basis. First assume that such a basis exists, i.e. $\sigma_1(\theta),\dots,\sigma_d(\theta)$ form a $\zed$-basis for $\O_K$ for some $\theta \in \O_K$. Notice that
$$\Sigma(\sigma_1(\theta)),\dots,\Sigma(\sigma_d(\theta))$$
forms a basis for $\Lambda_K$. Now for each $1 \leq j \leq d$,
\begin{eqnarray*}
\Sigma_K(\sigma_j(\theta)) & = & \Sigma_K(\sigma^{d-j+1}(\theta)) = \rho (\Sigma_K(\sigma^{d-j}(\theta))) \\
& = & \rho^2 (\Sigma_K(\sigma^{d-j-1}(\theta))) = \dots = \rho^{d-j} (\Sigma_K(\sigma(\theta))) = \rho^{d-j+1} (\Sigma_K(\theta)),
\end{eqnarray*}
by recursive application of~\eqref{rho_sigma}. Therefore $\Lambda_K$ is spanned by the basis
$$\Sigma_K(\theta),\rho(\Sigma_K(\theta))\dots,\rho^{d-1}(\Sigma_K(\theta)),$$
and hence it is simple cyclic.

Next suppose $\Lambda_K$ is simple cyclic, then 
$$\Lambda_K = \spn_{\zed} \{ \bx, \rho(\bx),\dots,\rho^{d-1}(\bx) \}$$
for some $\bx \in \Lambda_K$. Let $\theta \in \O_K$ be such that $\Sigma_K(\theta) = \bx$, then
$$\Sigma_K(\theta),\rho(\Sigma_K(\theta))\dots,\rho^{d-1}(\Sigma_K(\theta))$$
is a basis for $\Lambda_K$, where for each $1 \leq j \leq d$,
$$\Sigma_K(\sigma_j(\theta)) = \rho^{d-j+1} (\Sigma_K(\theta)),$$
as we derived above. Therefore $\sigma_1(\theta),\dots,\sigma_d(\theta)$ form a $\zed$-basis for $\O_K$, i.e. $K$ has a normal integral basis.
\endproof

An example of a family of such number fields comes from the cyclotomic fields.

\begin{cor} \label{cyclotomic} Let $K = \que(\zeta_n)$ be $n$-th cyclotomic field. There exists an ordering of the embeddings
$$\sigma_1,\dots,\sigma_{\phi(n)} : K \hookrightarrow \cee,$$
for which the lattice $\Lambda_K = \Sigma_K(\O_K)$ is cyclic if and only if $n=2, 4, p^k$ or $2p^k$ for some odd prime $p$ and positive integer $k$. Further, the cyclotomic lattice $\Lambda_K$ is simple cyclic if and only if $n = 2, p$, or $2p$ for an odd prime $p$.
\end{cor}

\proof
Recall that $K/\que$ is Galois with the Galois group $G \cong (\zed/n\zed)^{\times}$, which is cyclic if and only if $n=2, 4, p^k$ or $2p^k$ for some odd prime $p$ and positive integer $k$: these are precisely the values of $n$ for which primitive roots modulo~$n$ exist. If $G = \left< \sigma \right>$ is cyclic, order the embeddings as in~\eqref{order}, and the first statement follows from Lemma~\ref{nf_cyclic}. Further, the cyclotomic field $K = \que(\zeta_n)$ is tamely ramified if and only if $n$ is squarefree. Thus for $K$ to be cyclic {\it and} tamely ramified, $n$ must be equal to $2, p$, or $2p$ for an odd prime $p$. The second assertion then follows by Lemma~\ref{K_sc}.
\endproof

\proof[Proof of Theorem~\ref{Galois}] The theorem now follows from Lemmas~\ref{nf_cyclic}, \ref{K_sc} and Corollary~\ref{cyclotomic}.
\endproof

\begin{rem} \label{more_cyclic} It is known that the lattice $\Lambda_K$ is well-rounded if and only if $K$ is a cyclotomic field: this was proved in \cite{lf:petersen} for a slightly different embedding, but the argument is identical for our embedding $\Sigma_K$ (see Lemma~7.1.5 of~\cite{kitaoka}). Hence the only cyclic well-rounded lattices of the form $\Lambda_K$ are those characterized in Corollary~\ref{cyclotomic}, while there are also other cyclic lattices $\Lambda_K$ (from non-cyclotomic cyclic extensions) that are not well-rounded. Indeed, consider for example the real quadratic (hence cyclic) extension $\que(\sqrt{5})/\que$. It is tamely ramified, and hence has a normal integral basis
$$\frac{1+\sqrt{5}}{2}, \frac{1-\sqrt{5}}{2}.$$
Thus $\Lambda_K$ is simple cyclic by Lemma~\ref{K_sc}, however it is not well-rounded: its only minimal vectors are $\pm \begin{pmatrix} 1 \\ 1 \end{pmatrix}$ (on the other hand, it is stable). In fact, any quadratic number field $\que(\sqrt{D})$ for squarefree $D$ is cyclic, thus (when $D \neq -1,-3$) gives rise to a non-well-rounded cyclic lattice (although not necessarily simple cyclic like in the example above, e.g. if $D \not\equiv 1\ (\md 4)$ then $\que(\sqrt{D})$ is not tamely ramified).
\end{rem}

\bigskip

\noindent
{\bf Acknowledgement:} We wish to thank Pavel Guerzhoy for some helpful remarks. We are also very grateful to the anonymous referee for a thorough reading and multiple suggestions and corrections that improved the quality of this paper.

\bibliographystyle{plain}  

\end{document}